\documentclass{article}
\usepackage{graphicx,subcaption}  % Required for inserting images

\usepackage{xcolor}
\usepackage{amsfonts,amssymb,amsthm,amsmath}
\usepackage[margin=1in]{geometry}
\usepackage{tikz}
\usepackage{pgfplots}
\pgfplotsset{compat=1.8}

\definecolor{citegreen}{HTML}{208054}
\definecolor{citeblue}{HTML}{0055cc}

\usepackage{hyperref}
\hypersetup{colorlinks=true,citecolor=citeblue,linkcolor=citeblue,filecolor=citeblue,urlcolor=magenta,breaklinks=true}

\newcommand{\ZZ}{\mathbb{Z}}
\newcommand{\RR}{\mathbb{R}}
\newcommand{\CC}{\mathbb{C}}

\newcommand{\cM}{\mathcal{M}}

\newtheorem{theorem}{Theorem}
\newtheorem{proposition}{Proposition}

\newtheorem{lemma}{Lemma}

\theoremstyle{definition}
\newtheorem{definition}{Definition}

\renewcommand{\i}{\mathrm{i}}

\newcommand{\conv}{\text{conv}}

\newcommand{\Mon}{\mathcal{M}}

\newcommand{\kMon}[1]{\mathcal{M}(#1)}

\title{Exact algorithms for \\
quadratic optimization over roots of unity}
\date{\textsuperscript{\textdagger} Department of Applied Mathematics and Theoretical Physics, University of Cambridge,
United Kingdom }

\author{Ahmad Al-Sulami\textsuperscript{\textdagger} \and Hamza Fawzi\textsuperscript{\textdagger} \and Shengding Sun\textsuperscript{\textdagger}}

\begin{document}

\maketitle
\begin{abstract}
We consider the problem of optimizing a multivariate quadratic function where each decision variable is constrained to be a complex $m$'th root of unity. Such problems have applications in signal processing, MIMO detection, and the computation of ground states in statistical physics, among others.
Our contributions in this paper are twofold.
We first study the convergence of the sum-of-squares hierarchy and prove its convergence to the exact solution after only $\lfloor n/2\rfloor+1$ levels (as opposed to $n$ levels). Our proof follows and generalizes the techniques and results used for the binary $m=2$ case developed by Fawzi, Saunderson, Parrilo.  
Second, we construct an integer \emph{binary} reformulation of the problem based on zonotopes which reduces by half the number of binary variables in the simple reformulation. We show on numerical experiments that this reformulation can result in significant speedups (up to 10x) in solution time.
\end{abstract}
\section{Introduction}

Quadratic unconstrained binary optimization, henceforth QUBO, is the problem of maximizing a quadratic objective over binary variables.
This problem is most commonly stated in the following form:
\begin{equation} \label{eq:QUBO}
    \begin{aligned}
        \text{max } & x^T Qx\\
        \text{s.t. }&  x\in \{-1,1\}^n\\
    \end{aligned}
\end{equation}
where $Q$ is a symmetric matrix of size $n \times n$.
This NP-hard problem has attracted interest within the mathematical programming community, in part because it can model graph partitioning problems, but also because it has various applications in electrical engineering, statistical physics, and machine learning \cite{poljak1995maximum,ding2001min,barahona1988application}. There has been renewed interest in such problems recently, driven by the appearance of new quantum devices which are particularly tailored for problems of the form \eqref{eq:QUBO}. This has driven optimization researchers to find new performant classical algorithms for solving QUBO problems, see e.g., \cite{dunning2018works,madam,junger2021quantum,rehfeldt2023faster}. 

In this paper we are interested in the generalization of QUBO where the binary set $\{-1,1\}^n$ is replaced by the set of $m$'th roots of unity, namely:
\begin{equation} \label{eq:gen-QUBO}
    \begin{aligned}
        \text{max } & z^* Qz = \sum_{ij} Q_{ij} z_i^* z_j\\
        \text{s.t. }&  z\in B_m^n\\
    \end{aligned}
\end{equation}
where
\begin{equation}
    B_m^n := \{z \in \mathbb{C}^n: z_i^m = 1 \text{ for } i = 1,\dots,n\}.
\end{equation}
Here, $Q$ is a general $n\times n$ Hermitian matrix. Problems of the form \eqref{eq:gen-QUBO} appear in modern signal processing and maximal likelihood estimation problems \cite{so2008probabilistic}, including in angular synchronization \cite{bandeira2017tightness}, the multiple-input multiple-output detection problem \cite{jiang2021tightness,lu2019tightness,mobasher2007near}, phase retrieval \cite{waldspurger2015phase}, and unimodular radar code design~\cite{soltanalian2014designing,lu2020enhanced}. The problem \eqref{eq:gen-QUBO} is also directly related to the problem of computing the ground state of the well-known vector Potts model \cite{pottsmodel} from statistical physics, where the coefficients $Q_{ij}$ play the role of the (negative) coupling coefficients.

\subsection{Contributions}\label{sec:contributions}

In this paper, we focus on algorithms to solve \eqref{eq:gen-QUBO} exactly. Our contributions are twofold:

\paragraph{Tightness of Lasserre hierarchy}

Our first result concerns the Lasserre hierarchy to solve \eqref{eq:gen-QUBO}. The Lasserre hierarchy is a semidefinite-programming based hierarchy which gives increasingly accurate upper bounds on the solution of \eqref{eq:gen-QUBO}. We adopt in this paper the sum-of-squares point of view, wherein the $k$'th level of the hierarchy is defined as the best upper bound that can be obtained via degree $k$ sum-of-squares certificates. In order to introduce the hierarchy more precisely, we need a few technical definitions. First, note that any function $f:B_m^n \to \RR$ can be expanded in the monomial (Fourier) basis given by
\begin{equation*}
     \Mon^n_m:= \{z^\alpha = \prod_{i=1}^n z_i^{\alpha_i}:  \alpha \in \ZZ^n_m = \{0,\dots,m-1\}^n\}.
\end{equation*}
(In the binary case, $\Mon^n_m$ is the set of square-free monomials.) Due to the cyclic nature of $\ZZ_m$ the notion of degree of monomials is not clearly defined. Hence, we need to define the \emph{signature} for a monomial $z^\alpha$ which will allow us to grade the monomials in $\Mon^n_m$ and index the levels of hierarchy.
\begin{definition}\label{def:signature-intro}
    Given a monomial $z^\alpha\in \Mon^n_m$, its \emph{signature} is the ordered $(m-1)$ tuple $(n_1(\alpha),\ldots,n_{m-1}(\alpha))$, where $n_i(\alpha)=|\{j \in \{1,\ldots,n\}:\alpha_j=i\}|$ is the number of variables in $z^\alpha$ whose exponent is $i$. The set of all monomials with signature $(n_1,\ldots,n_{m-1})$ is denoted by $M(n_1,\ldots,n_{m-1})$. The \emph{level} $k$ monomials $\Mon^n_m(k)$ (or $\kMon{k}$ when $n$ and $m$ are clear from context) consist of monomials $z^\alpha$ such that $n_i(\alpha)\le k$ for all $1\le i\le m-1$. 
\end{definition}
It is obvious to check that in the case $m=2$, the signature of a monomial $z^\alpha$ is a single integer equal to its degree (as square-free monomial).
We are now ready to state the main theorem concerning the sum of squares relaxation:
\begin{theorem}\label{thm:main-intro}
Let $Q$ be an $n\times n$ Hermitian matrix and consider the function $f:B_m^n \to \RR$ defined by
\begin{equation}
    \label{eq:quadf}
    f(z) = z^* Q z = \sum_{1\leq i,j \leq n} Q_{ij} z_i^{-1} z_j.
\end{equation}
Assume that $f(z) \geq 0$ for all $z \in B_m^n$. Then $f$ admits a sum-of-squares certificate of the form $f(z)=\sum_{j=1}^r |g_j(z)|^2$ on $B_m^n$ where each $g_j:B_m^n\to \RR$ is supported on
\begin{equation}
    \label{eq:monoms-main-thm}
    \kMon{\lfloor n/2\rfloor +1} =\{z^\alpha:~ n_j(\alpha)\le \lfloor n/2\rfloor+1,\quad \forall 1\le j\le m-1\}.
\end{equation}
In short, the level-$k$ sum-of-squares hierarchy for \eqref{eq:gen-QUBO} is exact for $k = \lfloor n/2\rfloor +1$. 
\end{theorem}

As a comparison, \cite[Theorem 2]{fawzi2016sparse} proves that for $m=2$, the level $\lceil n/2\rceil$ sum-of-squares hierarchy is exact. Our result works for general $m$, and matches the same level bound for all odd $n$, and differs by one for all even $n$.
%\footnote{Given $f:B_m^n\to\RR,~f(z)=z^*Qz$ which is nonnegative on $B_m^n$, let $f_1:B_m^{n-1}\to \RR$ be obtained from $f$ by setting $z_n=1$. The main result we prove in Theorem \ref{thm:main} can be applied to $f_1$, showing that $f_1$ has a sum-of-squares certificate on $\mathcal{M}_m^{n-1}(k)$, where $k=\lfloor (n-1)/2\rfloor+1=\lceil n/2\rceil$, which would match the bound in \cite{fawzi2016sparse} for all $n$. However, we could not prove that any sum-of-squares certificate for $f_1$ supported on $\mathcal{M}_m^{n-1}(k)$ can be transformed into a sum-of-squares certificate for $f$ supported on $\mathcal{M}_m^n(k)$.}
While the techniques we use are similar to the previously cited paper and rely on the careful construction of chordal completions of the so-called Cayley graph, this generalization is highly nontrivial because the Cayley graph in our setting is significantly more complex than in the case where $m=2$.

We note that the result we are actually able to prove is considerably stronger than the statement in Theorem \ref{thm:main-intro}. What we show is that the level $\lfloor n/2\rfloor+1$ is exact for any polynomial $f$ which is a linear combination of monomials in at most two variables each, where degree can be arbitrary (see Theorem \ref{thm:main}). However we are not aware of optimization applications beyond the quadratic case so that is why we decided to focus only on functions of the form \eqref{eq:quadf}.
 We also note that a previous generalization of the theorem in \cite{fawzi2016sparse} was given in \cite{sakaue2017exact}, however the latter concerned nonquadratic functions $f$ in binary variables.

\paragraph{An exact real binary reformulation of \eqref{eq:gen-QUBO}} The second contribution of this paper concerns the practical computation of the optimal value of problems of the form \eqref{eq:gen-QUBO}. Noting that most existing integer programming solvers do not natively support complex discrete variables in $B_m$, but many do support quadratic objectives for binary variables, (e.g., Gurobi \cite{gurobi}, SCIP \cite{bolusani2024scip}, CPLEX \cite{cplex2009v12}, etc.) we study the question of what the best way is to reformulate the problem \eqref{eq:gen-QUBO} as a purely binary quadratic optimization problem, with possibly additional linear constraints. We first present a simple reformulation of \eqref{eq:gen-QUBO} as a binary maximization problem in $mn$ variables and $n$ linear equality constraints. We then proceed to construct a nontrivial reformulation in the case of $m$ even, which uses only $mn/2$ binary variables.
Our construction is based on the observation that the regular $m$-gon, in the case $m$ even is a \emph{zonotope}, i.e., it is the linear projection of a cube in higher dimensions. Our construction can be summarized in the following statement.
\begin{theorem}[Informal, see Theorem \ref{thm:zonotope-reformulation} for detailed statement]

For any integer $n$ and even $m$, there exists an exact reformulation of \eqref{eq:gen-QUBO} as a pure binary quadratic problem with $nm/2$ binary variables and with $n$ convex inequalities.
\end{theorem}
We have tested our new reformulation on problems from MIMO detection and statistical physics, and we have observed that the new reformulation can achieve significant speedups (up to 10x) compared to the simple reformulation.

\subsection{Related work}

The problem \eqref{eq:gen-QUBO} is sometimes known as Discrete Complex Quadratic Optimization \cite{zhang2006complex,huang2010approximation,so2007approximating}, 
and there has been a lot of work on semidefinite optimization-based approximation algorithms for \eqref{eq:gen-QUBO}. The starting point of these relaxations is the celebrated Goemans-Williamson relaxation in the case $m=2$ which is defined by
\begin{equation}
\label{eq:QUBO-SDP}
\begin{aligned}
        \text{max } & \text{tr}(QX)\\
        \text{s.t. }&  X_{ii}=1,\quad 1\le i\le n\\
        & X\succeq 0,
\end{aligned}
\end{equation}
where $X$ is real symmetric.
With randomized rounding, approximation guarantees on this relaxation have been obtained for different classes of objective matrices $Q$ in the case $m=2$, see for example \cite[Chapter 2]{blekherman2012semidefinite}.

The SDP \eqref{eq:QUBO-SDP} for any $m\ge 3$ can be strengthened by adding the additional constraints $X_{ij}\in \conv(B_m)$ for all $i,j$ (where now $X$ is assumed Hermitian). Goemans and Williamson showed in \cite{goemans2001approximation} that for $m=3$ and Max-3-Cut objectives this strengthened SDP has approximation factor $\approx 0.836$.
Different rounding approaches for larger values of $m \geq 3$ have also been proposed, see \cite{zhang2006complex,huang2010approximation} and \cite{so2008probabilistic,jiang2021tightness,lu2019tightness} in particular for the MIMO detection problem. \\
Recently, Sinjorgo et al. in \cite{sinjorgo2024cuts} studied problem \eqref{eq:gen-QUBO} from a geometric point of view, and obtained valid linear inequalities on the polytope $\conv \left\{ zz^* : z \in B_m^n \right\}$ that can be used to tighten \eqref{eq:QUBO-SDP} for general $m$.
The Lasserre/sum-of-squares hierarchy for the $m=2$ QUBO case \eqref{eq:QUBO} has been widely studied \cite{laurent2004semidefinite,fawzi2016sparse,slot2023sum}. The $m\to \infty$ limit coincides with trigonometric sum-of squares hierarchy, which has been studied in \cite{naftalevich2006trigonometric} and recently in \cite{bach2023exponential}. 
To the best of our knowledge, there has not been work that explicitly studies the sum-of-squares hierarchy for general finite $m$.

\subsection{Paper organization}
This paper is organized as follows. We provide the necessary preliminary background in Section \ref{sec:preliminaries}, consisting of a discussion on chordal graphs, and another on zonotopes. In Section \ref{sec:m-sos}, we prove our main result on the convergence of the Lasserre hierarchy (Theorem \ref{thm:main}). In Section \ref{sec:reformulations}, we discuss the reformulation of \eqref{eq:gen-QUBO} as a real binary program, including a standard reformulation for general $m$, and a more efficient reformulation for even $m$ based on zonotopes (Theorem \ref{thm:zonotope-reformulation}). Section \ref{sec:numerics} presents numerical examples for the reformulations in Section  \ref{sec:reformulations}.
Section \ref{sec:conclusions} concludes the paper and discusses some related open problems. 

\subsection{Notations}
For $z \in \CC$, we denote by $z^*$ its complex conjugate. By extension, for a complex vector $z = [z_1,\dots, z_n]\in \CC^n$, we denote its conjugate transpose by $z^* = [z_1^*,\dots,z_n^*]^T$, and for a complex matrix $ Z \in \CC^{n\times n}$, we denote its Hermitian conjugate by $Z^*$. A matrix $Z \in \CC^{n\times n}$ is called Hermitian if $Z^* = Z$. The set of all Hermitian matrices of size $n\times n$ is denoted by $\mathbb{H}^n$, and its subset of positive semidefinite matrices forms a cone that we denote by $\mathbb{H}^n_+$. In the case of real matrices $X \in \RR^{n \times n}$, we say that it is symmetric if $X^T = X$, and we denote the set of symmetric matrices by $\mathbb{S}^n$ and the cone of symmetric positive semidefinite matrices by $\mathbb{S}^n_+$. We use the notation $X \succeq 0$ to denote that $X$ is symmetric positive semidefinite.  We adopt the multi-index notation for monomials: given $z \in \CC^n$ and $\alpha \in \ZZ^n$, we define the monomial $z^\alpha = \prod_{i=1}^n z_i^{\alpha_i}$, and also define $|\alpha| = \sum_{i=1}^n |\alpha_i|$. We denote the floor of a real number $x$ by $\lfloor x\rfloor$, and its ceiling by $\lceil x\rceil$.\\

\section{Preliminaries}\label{sec:preliminaries}
In this section, our aim is to cover the background necessary to follow the rest of the paper.  

\subsection{Chordal graphs and the running intersection property}
Let $G = (V,E)$ be an undirected graph. We say that $G$ is \emph{chordal} if any cycle in $G$ of length greater than three has a chord. We say that $G'$ is a \emph{chordal cover} of $G$ if $G' = (V,E')$, with $E \subseteq E'$, is chordal. A \emph{clique} $C \subseteq V$ of $G$ is a subset of the vertices such that every two distinct vertices in $C$ are adjacent in $G$. In other words, $C$ is a clique of $G$ if $\{i,j\} \in E$ for all $i,j \in C$ such that $i \neq j$. If $C$ is not a strict subset of another clique $C'$ of $G$, then we call it a \emph{maximal clique}. The clique sometimes also refers to its induced subgraph, consisting of the vertices in the clique and all edges between these vertices.

We are mainly interested in the above graph theoretical concepts because they are useful in the context of semidefinite programming. It is known, see e.g., \cite[Theorem 2.3]{agler1988positive} that if $Q$ is a positive semidefinite matrix which has a chordal sparsity pattern $G=(V,E)$, then $Q$ can be decomposed as a sum of positive semidefinite matrices which are supported on the maximal cliques of $G$.

There are many different ways to characterize chordal graphs. The following lemma which shows how chordal graphs can be obtained by clique sum operations (gluing together two subgraphs on a clique) starting from cliques will be especially useful for us. It is known as the \emph{running intersection property}, see e.g., \cite{lasserre2006convergent}. 

\begin{lemma}\label{lem:running_intersection}
    Let $V_0,V_1,\ldots,V_k$ be subsets of the vertex set $V$, so that $V=\bigcup_{i=0}^k V_i$ and $V_0,\ldots,V_k$ satisfy the running intersection property: for any $1\le j\le k$, there exists some $s<j$ such that
    \[
    V_j\cap\left(\bigcup_{i=0}^{j-1}V_i\right)\subseteq V_s. 
    \]
    Then the graph with vertex set $V$ and edge set $\bigcup_{i=0}^k \mathcal{C}(V_i)$ is chordal, where $\mathcal{C}(V_i)$ is the clique with vertex set $V_i$. 
\end{lemma}

\subsection{Zonotopes}

A \emph{zonotope} $Z$ in $\RR^m$ is a polytope that can be decomposed as the Minkowski sum of line segments, that is: 
\begin{equation*}
    Z = L_1+\dots+L_d , \hspace{5pt} L_i = [a_i,b_i] \text{ with } a_i,b_i \in \RR^m, \text{ for } i=1,\dots,d.
\end{equation*}
Alternatively, a zonotope is a polytope that can be written as the affine image of a cube $[0,1]^d$ for some $m$. The affine map corresponding to the representation above is the map $T(\lambda) = \sum_{i=1}^d a_i + \lambda_i (b_i - a_i)$.
% A simple example of a zonotope would be the convex hull of the hypercube:
% \begin{equation*}
%     \conv(\{-1,1\}^n) = \sum_{i=1}^n [-e_i,e_i],
% \end{equation*}
% where $e_i$ is the $i$-th canonical basis vector in $\RR^n$ with the entry $1$ in the $i$-th position, and zeros elsewhere.
Later, we will see how to construct a lift for $\conv(B_m^n)$ with even $m$ using zonotopes. 

\section{Moment-sum-of-squares hierarchy}\label{sec:m-sos}

In this section we study the moment/sum-of-squares (moment-SOS) hierarchy for \eqref{eq:gen-QUBO}. The moment-SOS hierarchy is a powerful tool to obtain increasingly tight semidefinite relaxations for polynomial and discrete optimization problems. It can be defined in different ways, notably either from the moment/primal point of view, or from sum-of-squares/dual point of view. In this work we focus on the sum-of-squares point of view. Then, the $k$'th level of the moment-SOS hierarchy \eqref{eq:gen-QUBO} is defined as the solution of the following optimization problem:
\begin{equation}
\label{eq:ksosproblem}
\begin{array}{ll}
\min & t\\
\text{s.t.} & f(z) = t - z^*Qz \text{ is a $k$-sum-of-squares}
\end{array}
\end{equation}
where a function $f(z)$ is called a $k$-sum-of-squares if it can be expressed as $f(z) = \sum_{j} |g_j(z)|^2$ for some functions $g_j:B_m^n \to \CC$ which are linear combinations of monomials from $\cM(k)$ (recall the definition of $\cM(k)$ from Definition \ref{def:signature-intro}). It is by now a standard fact that \eqref{eq:ksosproblem} can be expressed as a semidefinite program with a positive semidefinite constraint of size $|\cM(k)|\times |\cM(k)|$. We omit a proof of this fact in this manuscript and refer e.g., to \cite{fawzi2016sparse}.

For binary quadratic programming, \cite{fawzi2016sparse} shows that the level $\lceil n/2\rceil$ relaxation is guaranteed to be exact. In fact, the paper \cite{fawzi2016sparse} proposes a general framework to analyze Fourier sum-of-squares on finite abelian groups. Given a finite abelian group $G$, let $\hat{G}$ denote its characters (dual group) so that any function $f:G\to\CC$ has a Fourier decomposition $f(z)=\sum_{\chi\in \hat{G}}\hat{f}(\chi)\chi(z)$ for all $z \in G$. The main result of \cite{fawzi2016sparse} is stated as follows: 

\begin{theorem}[Theorem 1, \cite{fawzi2016sparse}]\label{thm:fsos}
    Let $S\subseteq \hat{G}$, let $\Gamma$ be a chordal cover of $\text{Cay}(\hat{G},S)$, and for each maximal clique $\mathcal{C}$ of $\Gamma$, let $\chi_{\mathcal{C}}$ be an element of $\hat{G}$. Define
    \[
    \mathcal{T}:=\mathcal{T}(\Gamma,\{\chi_{\mathcal{C}}\})=\bigcup_{\mathcal{C}} \chi_{\mathcal{C}} \mathcal{C},
    \]
    where the union is over all the maximal cliques of $\Gamma$ and where $\chi_{\mathcal{C}} \mathcal{C}:=\{\chi_{\mathcal{C}}\chi:~\chi\in \mathcal{C}\}$ is the translation of $\mathcal{C}$ by $\chi_{\mathcal{C}}$. Then any nonnegative function with support $S$ admits a sum-of-squares certificate with support $\mathcal{T}$. 
\end{theorem}

We apply this result to the case $G=B_m^n\simeq \ZZ_m^n$. Its dual group is given by the monomials:
$$\hat{G} = \Mon_m^n =\{z^\alpha: \alpha \in \{0,1,\dots,m-1\}^n\},$$ 
which are group homomorphisms on $\ZZ_m^n$. \\
Note that the usual notion of degree of monomials does not work very well with $B_m^n$. Instead, we grade the monomials by their signatures which were already defined in Definition \ref{def:signature-intro}. For reader's convenience, recall that the signature of a monomial $z^\alpha$ is the ordered $(m-1)$ tuple $(n_1(\alpha),\ldots,n_{m-1}(\alpha))$, where $n_i(\alpha)=|\{j:\alpha_j=i\}|$ is the number of variables in $z^\alpha$ with exponent $i$. The set of all monomials with signature $(n_1,\ldots,n_{m-1})$ is denoted by $M(n_1,\ldots,n_{m-1})$. The \emph{level} $k$ monomials $\kMon{k}$ consist of monomials $z^\alpha$ such that $n_i(\alpha)\le k$ for all $1\le i\le m-1$. 

For $m=2$ the level $k$ monomials are precisely ones with degree at most $k$, when expressed as square-free monomials. For $m=3$, as $z_j^2=z_j^{-1}$ for every $j$, the elements in $\Mon_m^n$ can be written as $z^\alpha z^{-\beta}$ where $\alpha,\beta\in \{0,1\}^n$ have distinct entries of ones. Then the level $k$ monomials can be expressed by $\{z^\alpha z^{-\beta}:~|\alpha|\le k,|\beta|\le k\}$. This agrees with the truncation rule proposed in \cite{josz2018lasserre} for general complex polynomial optimization.

Now we turn to the proof of Theorem \ref{thm:main-intro}. As announced in section \ref{sec:contributions}, we will prove a considerably stronger statement for linear combinations of monomials in at most two variables each. 

\begin{theorem}\label{thm:main}
    Let $f:B_m^n\to\RR$ be a nonnegative polynomial such that 
    \[
    f\in\text{span}(\{z^\alpha:~ |\{i:\alpha_i\ne 0\}|\le 2\}).
    \]
    Then $f=\sum_{j=1}^r g_j \bar{g}_j$ on $B_m^n$ where each $g_j$ is supported on
    \[
    \kMon{\lfloor n/2\rfloor+1}:=\{z^\alpha:~ n_j(\alpha)\le \lfloor n/2\rfloor+1,\quad \forall 1\le j\le m-1\}.
    \]
\end{theorem}

\begin{proof}
    We use Theorem \ref{thm:fsos} with $G=B_m^n$ and $S=\{z^\alpha:~ |\{i:\alpha_i\ne 0\}|\le 2\}$. Recall $n_j(\beta)=|\{i:~\beta_i=j\}|$. We first prove the following lemma which will be used later in our proof. 

    \begin{lemma}\label{lem:Cay-adjacent}
        For any $\beta,\gamma\in\{0,\ldots,m-1\}^n$, if $z^\beta$ and $z^\gamma$ are adjacent on $\text{Cay}(\hat{G},S)$, then $|n_j(\gamma)-n_j(\beta)|\le 2$ for all $1\le j\le m-1$. 
    \end{lemma}
    \begin{proof}[Proof of Lemma \ref{lem:Cay-adjacent}]
    By definition of $\text{Cay}(\hat{G},S)$, $z^\beta$ and $z^\gamma$ are adjacent on $\text{Cay}(\hat{G},S)$ if and only if $\beta$ and $\gamma$ differ in at most two coordinates, i.e., there exist $i,i'\in \{1,\ldots,n\}$ and $c_1,c_2\in \{0,\ldots,m-1\}$ such that $\gamma=\beta+c_1e_i+c_2e_{i'}$, where $e_i$ is the $i$-th indicator vector which has one at coordinate $i$ and zero everywhere else. 

    For any $\beta\in \{0,\ldots,m-1\}^n$, we first show $|n_j(ce_i+\beta)-n_j(\beta)|\le 1$ for all $1\le j\le m-1$, where $i\in \{1,\ldots,n\},c\in \{0,\ldots,m-1\}$. When $c = 0$, the statment is trivial as $n_j(ce_i+\beta)-n_j(\beta) = 0$ for all $1\le j\le m-1$. When $c\ne 0$ we have

    \begin{equation*}
        n_j(ce_i+\beta)-n_j(\beta)=\begin{cases}
            1 & j\equiv \beta_i+c\mod m,\\
            -1 & j=\beta_i,\\
            0 & \text{otherwise.}
        \end{cases}
    \end{equation*}
    Now to show $|n_j(\gamma)-n_j(\beta)|\le 2$ for all $1\le j\le m-1$ where $\gamma=\beta+c_1e_i+c_2e_{i'}$, we use the above statement twice and triangle inequality. We have
    \begin{equation*}
        \begin{aligned}
            |n_j(c_1e_i+c_2e_{i'}+\beta)-n_j(\beta)| & = |n_j(c_1e_i+c_2e_{i'}+\beta)-n_j(c_1e_i+\beta) + n_j(c_1e_i+\beta)-n_j(\beta)|\\
            & \le |n_j(c_1e_i+c_2e_{i'}+\beta)-n_j(c_1e_i+\beta)|+|n_j(c_1e_i+\beta)-n_j(\beta)| \\
            & \le 2.
        \end{aligned}
    \end{equation*}

    The lemma is thus proved. 
    \end{proof}

     Now let $k=\lfloor n/2\rfloor +1$, which is the smallest integer $k$ such that $2k\ge n+1$. We define the following subsets of monomials $\{V_0,\ldots,V_{m-1}\}$: 
\begin{equation}
\label{eq:defVj}
\begin{aligned}
V_0 &:=\{z^\alpha:~ n_i(\alpha)\le k\text{ for all }1\le i\le m-1\},\\
V_j&:=\{z^\alpha:~ n_j(\alpha)\ge k-1\},\quad 1\le j\le m-1.
\end{aligned}
\end{equation}
A Venn diagram to illustrate these sets in the case $m=5$ is given in Figure \ref{fig:venn_diagram_Vs}. We clearly have $\bigcup_{j=0}^{m-1} V_j = \hat{G}$ (note the $V_j$'s are not disjoint). Let $\Gamma$ be the graph with vertex set $\hat{G}$ and edge set equal to $\bigcup_{j=0}^{m-1} \mathcal{C}(V_j)$ where $\mathcal{C}(V_j)$ is the clique with vertex set $V_j$
\[
\Gamma = ( \hat{G} , \bigcup_{j=0}^{m-1} \mathcal{C}(V_j) ).
\]
We first prove that all the edges of $\text{Cay}(\hat{G},S)$ are covered by $\Gamma$.
\begin{proposition}\label{prop:cover}

$\text{Cay}(\hat{G},S)$ is a subgraph of $\Gamma$.
\end{proposition}
\begin{proof}[Proof of Proposition \ref{prop:cover}]
    Let $\{z^\alpha,z^\beta\}$ be an edge in $\text{Cay}(\hat{G},S)$. We will show that the edge $\{z^\alpha,z^\beta\}$ is included in one of $\mathcal{C}(V_j)$ for some $0\leq j\leq m-1$. Since $\bigcup_{j=0}^{m-1}V_i=\hat{G}$, there exists $p,q\in \{0,\ldots,m-1\}$ such that $z^\alpha\in V_p,z^\beta\in V_q$. If $p=q$ then the edge $\{z^\alpha,z^\beta\}$ is covered by $\mathcal{C}(V_p)$. The same is true if $z^\alpha \in V_p \cap V_q$ or $z^\beta \in V_p \cap V_q$. If not, then this means that $p \neq q$ and $z^\alpha\in V_p\setminus V_q,z^\beta\in V_q\setminus V_p$. We proceed as follows: 
    \begin{itemize}
        \item We first show that necessarily $p\neq 0$ and $q\neq 0$. Indeed, assume for contradiction that $q=0$ (without loss of generality). Then we will show that $|n_p(\beta) - n_p(\alpha)| \geq 3$ which means that $\{z^\alpha,z^\beta\}$ cannot be an edge in $\text{Cay}(\hat{G},S)$. Note $z^\alpha\in V_0\setminus V_p$ implies $n_p(\alpha)\le k-2$, and $z^\beta\in V_p\setminus V_0$ implies $n_p(\beta)\ge k-1$ and there exists some $j\in\{1,\ldots,m-1\}$ such that $n_j(\beta)\ge k+1$. This implies $j=p$, since otherwise we would have $n = \sum_{i=0}^{m-1} n_i(\beta) \ge n_j(\beta)+n_p(\beta)\ge 2k>n$, a contradiction. Thus we get $n_p(\beta)\ge k+1$, and therefore $|n_p(\beta)-n_p(\alpha)|\ge 3$. This shows that $z^\alpha$ and $z^\beta$ are not adjacent in $\text{Cay}(\hat{G},S)$ by Lemma \ref{lem:Cay-adjacent}. 
        
        \item From the above, it means that $p\neq 0$ and $q \neq 0$. We will now show that the edge $\{\alpha,\beta\}$ in $\text{Cay}(\hat{G},S)$ is necessarily in $\mathcal{C}(V_0)$, in other words $z^{\alpha},z^{\beta} \in V_0$. From $z^{\alpha} \in V_p \setminus V_q$ and $z^{\beta} \in V_q \setminus V_p$, we get 
        \[
        n_p(\alpha)\ge k-1,\; n_q(\alpha)\le k-2, \; n_q(\beta)\ge k-1, \; n_p(\beta)\le k-2.
        \]
        Since $n_p(\alpha)\ge k-1$, for any $j\ne p$ we have $n_j(\alpha)\le n-n_p(\alpha)\le n-k+1\le k$. Similarly $n_q(\beta)\ge k-1$ and $n_j(\beta)\le k$ for any $j\ne q$. To prove $z^\alpha,z^\beta\in V_0$, it remains to show $n_p(\alpha)\le k$ and $n_q(\beta)\le k$. As $z^\alpha$ and $z^\beta$ are adjacent on $\text{Cay}(\hat{G},S)$, from Lemma \ref{lem:Cay-adjacent} we have $|n_j(\alpha)-n_j(\beta)|\le 2$ for all $j\in\{1,\ldots,m-1\}$, and hence $n_p(\alpha)\le n_p(\beta)+2\le k,n_q(\beta)\le n_q(\alpha)+2\le k$. 
    \end{itemize}
\end{proof}
\begin{figure}[ht]
    \centering
    \includegraphics[width=0.3\linewidth]{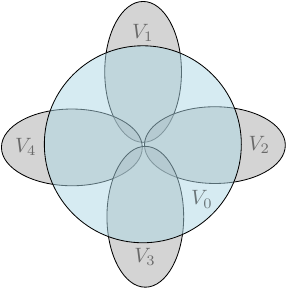}
    \caption{A simplified Venn diagram for the sets $V_0,V_1,\dots, V_4$ ($m = 5$)}
    \label{fig:venn_diagram_Vs}
\end{figure}
We now prove that the graph $\Gamma$ we constructed is chordal. 

\begin{proposition}\label{prop:chordal}
    $\Gamma=(\hat{G},\bigcup_{j=0}^{m-1} \mathcal{C}(V_j))$ is chordal.
\end{proposition}
\begin{proof}[Proof of Proposition \ref{prop:chordal}]
    To show $\Gamma$ is chordal, note for any distinct $p,q\in \{1,\ldots,m-1\}$ we have $V_p\cap V_q=\{z^\alpha\in\hat{G}:n_p(\alpha)\ge k-1,n_q(\alpha)\ge k-1\}$. From $n_p(\alpha)\ge k-1$ we get that $n_j(\alpha)\le n-n_p(\alpha)\le n-k+1\le k$ for all $j\neq p$ and $n_q(\alpha)\ge k-1$ implies $n_p(\alpha)\le k$. Thus $V_p\cap V_q\subseteq V_0$. Therefore $V_j\cap\left(\bigcup_{i=0}^{j-1}V_i\right)\subseteq \bigcup_{i=0}^{j-1} (V_j \cap V_i) \subseteq V_0$ for all $1\le j\le m-1$, and $\{V_0,\ldots,V_{m-1}\}$ satisfies the running intersection property. 
\end{proof}
To complete the proof of Theorem \ref{thm:main}, we need to show that for each $0 \leq j \leq m-1$ there exists $\chi_j \in \hat{G}$ such that $\chi_j V_j \subseteq \cM(k) = V_0$. For $j=0$ we can just let $\chi_j=1$ be the degree zero monomial, i.e. the identity element in $\hat{G}$. For all $1\le j\le m-1$ let $\chi_j=z_1^{-j}\ldots z_n^{-j}$. Since $2k\ge n+1$ we have
    \[\chi_j V_j=\{z^\alpha:~ n_0(\alpha)\ge k-1\}=\{z^\alpha:~ \sum_{i=1}^{m-1}n_i(\alpha)\le n-k+1\}\subseteq V_0.\]
    Thus Theorem \ref{thm:fsos} holds with $\mathcal{T}=V_0=\kMon{k}$. 
\end{proof}

\section{Reformulation into real binary programming}\label{sec:reformulations}

In this section, we discuss two ways to reformulate problem \eqref{eq:gen-QUBO} into real binary programming instances. We first discuss a simple reformulation that expresses problem \eqref{eq:gen-QUBO} as a binary quadratic problem with $mn$ binary variables. Then, we show how one can improve this construction to halve the number of binary variables in the case where $m$ is even. In the following section, we validate the usefulness of this new reformulation on  numerical experiments.

\subsection{A reformulation with \texorpdfstring{$mn$}{mn} variables}

Let $Q$ be a Hermitian $n\times n$ matrix and consider the problem
\begin{equation} \label{eq:gen-QUBO-3}
    \begin{aligned}
        \text{max } & z^* Qz = \sum_{ij} Q_{ij} z_i^* z_j\\
        \text{s.t. }&  z\in B_m^n.\\
    \end{aligned}
\end{equation}
We seek a reformulation of this discrete optimization that uses only binary variables. Let $\omega=\exp(2\pi\i/m)$ be the primitive $m$-th root. For each $j=1,\ldots,n$ and $k=0,\ldots,m-1$, we define binary variables $x_{j,k}=1$ if $z_j=\omega^k$ and zero otherwise. Then note that for any $j,l$:
\[
z^{-1}_jz_l=\begin{pmatrix}
    x_{j,0} & \ldots &x_{j,m-1}
\end{pmatrix}\begin{bmatrix}
    1 & \omega & \ldots &\omega^{m-1}\\
    \omega^{-1} & 1 & \ldots & \omega^{m-2}\\
    \vdots & \vdots & \ddots & \vdots\\
    \omega^{-(m-1)} & \omega^{-(m-2)} & \ldots & 1
\end{bmatrix}\begin{pmatrix}
    x_{l,0} \\ \vdots \\ x_{l,m-1}
\end{pmatrix}.
\]

Let $W$ be the matrix $W_{jl}=\omega^{l-j}$ and encode all $x_{j,k}$ variables into a binary vector $x$ of length $mn$ as follows $x=(x_{1,0},\ldots,x_{1,m-1},\ldots,x_{n,0},\ldots,x_{n,m-1})^T$. Then we have the following real binary reformulation of \eqref{eq:gen-QUBO}:

\begin{equation}\label{eq:basic-reformulation}
\max \{ z^* Q z : z \in B_m^n \} \; = \;
    \begin{array}[t]{ll}
        \underset{x=(x_{j,k})}{\text{max}} & x^\top \text{Re}(Q\otimes W)x\\
        \text{s.t. }& \sum_{k=0}^{m-1} x_{j,k}=1,\quad \forall 1\le j\le n,\\
        & v\in \{0,1\}^{mn}.
    \end{array}
\end{equation}

The equality constraints in the formulation above ensure that there is exactly one variable out of $\{x_{j,0},\ldots,x_{j,m-1}\}$ which is equal to one, for each $1\leq j\leq n$.

\subsection{A better reformulation}

In this section we show how to obtain a reformulation of \eqref{eq:gen-QUBO} as a binary maximization problem with only $mn/2$ binary variables (as opposed to $mn$) in the case where $m$ is even.

We will first consider the case $m=4$ and observe that one can get a reformulation using simply $2n$ binary variables (and no extra constraints).
%, instead of the $4n$ binary variables from the reformulation \eqref{eq:basic-reformulation}.

\paragraph{Case $m=4$} When $m=4$ our problem becomes
\begin{equation*}
\begin{aligned}
    \text{max }  & z^*Qz\\
    \text{s.t. } & z\in \{-1,1,-\i,\i\}^n.
\end{aligned}
\end{equation*}
Consider a change of variables $w=(1+\i) z$. Then the problem can be exactly rewritten as
\begin{equation*}
\begin{aligned}
    \text{max } & \frac{1}{2}w^*Qw\\
    \text{s.t. } & w\in \{\pm 1\pm \i\}^n
\end{aligned}
\end{equation*}
Let $w=x+y\i$ where $x,y\in\{-1,1\}^n$. For $v=\begin{pmatrix}
    x\\y
\end{pmatrix}$, we thus obtain
\begin{equation}
\max \{ z^* Q z : z \in B_4^n \} \; = \;
\begin{array}[t]{ll}
    \text{max } & \frac{1}{2}v^\top \begin{bmatrix}
        \text{Re}(Q) & -\text{Im}(Q)\\
        \text{Im}(Q) & \text{Re}(Q)
    \end{bmatrix}v\\
    \text{s.t. } & v\in \{-1,1\}^{2n}
\end{array}
\end{equation}
Thus we see that when $m=4$ the problem \eqref{eq:gen-QUBO} can be reformulated as a pure QUBO with only $2n$ variables.

\paragraph{General case for even $m$}

Now we show that for even $m$, one can get a reformulation of \eqref{eq:gen-QUBO} with only $mn/2$ binary variables and $n$ convex inequality constraints. Our main result is the following.

\begin{theorem}[Reformulation for even $m$]
\label{thm:zonotope-reformulation}
Let $m$ be even. For any Hermitian $n\times n$ matrix $Q$ we have
\begin{equation}\label{eq:zonotope-reformulation}
\max_{z \in B_m^n} z^* Q z = \max_{x \in \{-1,1\}^{nm/2}} \left\{ x^T \tilde{Q} x \text{ s.t. } \sum_{k=1}^{m/2-1} |x_{i(m/2)+k} - x_{i(m/2)+k+1}| \leq 2 ,\; \forall 0\leq i \leq n-1\right\}
\end{equation}
where $\tilde{Q} \in \RR^{nm/2 \times nm/2}$ is defined as $\tilde Q=(I_n\otimes A)^TQ_R(I_n\otimes A)$, where
\[
Q_R=\begin{bmatrix}
        Q_R(1,1) & \ldots & Q_R(1,n)\\
        \vdots & \ddots & \vdots\\
        Q_R(n,1) & \ldots & Q_R(n,n)
    \end{bmatrix}  \in \RR^{2n \times 2n} ,\quad  Q_R(i,j)=\begin{bmatrix}
        \text{Re}(Q_{ij}) & -\text{Im}(Q_{ij})\\
        \text{Im}(Q_{ij}) & \text{Re}(Q_{ij})
    \end{bmatrix} \in \RR^{2\times 2},
    \]
    \[
    A=\sin(\pi/m)\begin{bmatrix}
        \sin(\pi/m) & \sin(3\pi/m) & \ldots & \sin((m-1)\pi/m)\\
        -\cos(\pi/m) & -\cos(3\pi/m) & \ldots & -\cos((m-1)\pi/m)
    \end{bmatrix} \in \RR^{2\times m/2},
    \]
and
\[ I_n\otimes A=\begin{bmatrix}
        A & &\\& \ddots &\\& & A
    \end{bmatrix} \in \RR^{2n\times mn/2}. \]
\end{theorem}

The proof of this theorem relies on the fact that the regular $m$-gon, in the case $m$ is even, is a zonotope, i.e., it can be expressed as the projection of a cube in $m/2$ dimensions. More precisely, there is a linear map $T:\RR^{m/2}\to \CC$ such that $T([-1,1]^{m/2}) = \conv(B_m)$ where $B_m$ are the $m$'th roots of unity.

In order to prove our theorem, we need to know precisely which vertices from the cube in $[-1,1]^{m/2}$ map to the vertices of the regular $m$-gon. This is the object of the next lemma which is crucial for the proof of Theorem \ref{thm:zonotope-reformulation}.

\begin{lemma}\label{lem:zonotope_map}
    Let $m$ be an even integer, and let $B_m = \{z \in \CC : z^m=1\}$ be the $m$'th roots of unity. Then we have:
    \[
    B_m = \left\{T \epsilon : \epsilon \in \{-1,1\}^{m/2}:~\sum_{j=1}^{m/2-1} |\epsilon_j - \epsilon_{j+1}| \leq 2\right\}
    \]
    where $T:\CC^{m/2}\to \CC$ is the linear map described by
    \begin{equation}
        \label{eq:defTek}
    T(e_k) = -\i \sin(\pi/m) \exp\left[\frac{(2k-1)\pi\i}{m}\right],\quad 1\le k\le \frac{m}{2}.
    \end{equation}
\end{lemma}
\begin{proof}
    Note that the set of binary vectors $\epsilon \in \{-1,1\}^{m/2}$ that satisfy the total variation bound
    \[
    \sum_{j=1}^{m/2-1} |\epsilon_j - \epsilon_{j+1}| \leq 2
    \]
are precisely the vectors of the form $\pm w_k$ where
\[
w_k = (\underbrace{-1,\ldots,-1}_{k},\underbrace{1,\ldots,1}_{m/2-k}), \qquad 1\leq k \leq m/2.
\]
Note that we can also write $w_k$ as
\[
w_k = -\sum_{j=1}^k e_j + \sum_{j=k+1}^{m/2} e_j.
\]
Our goal is to show that $T(w_k) = \exp(2ik\pi /m)$ for all $1\leq k \leq m/2$. We start by writing
\begin{equation}
\label{eq:Twk}
T(w_k) = -\sum_{j=1}^k T(e_j) + \sum_{j=k+1}^{m/2} T(e_j).
\end{equation}
The following trigonometric identities are immediate to verify using geometric series: 
\[
\sum_{j=1}^k \exp\left[\frac{(2j-1)\pi\i}{m}\right] = \frac{e^{2\i k \pi /m} - 1}{2 \i \sin(\pi/m)}
\quad
\text{ and }
\quad
\sum_{j=k+1}^{m/2} \exp\left[\frac{(2j-1)\pi\i}{m}\right] = -\frac{e^{2\i k \pi /m} + 1}{2 \i \sin(\pi/m)}.
\]

Plugging in \eqref{eq:Twk}, using the definition of $T(e_k)$ from \eqref{eq:defTek}, we get $T(w_k) = \exp(2ik\pi/m)$ as desired.
\end{proof}

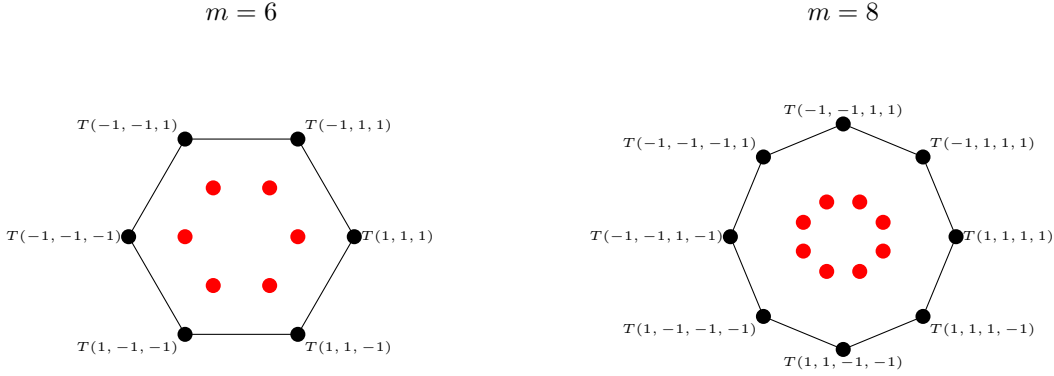
\begin{figure}[ht]
    \centering
\begin{tikzpicture}
    % m = 6
    \draw (0,3) node[] {$m=6$};
  \foreach \i in {0,1,...,5} {
    % Vertex hexagon
    \coordinate (h\i) at (60*\i:1.5); % radius = 1.5
    \coordinate (hi\i) at (60*\i:0.75); % radius = 0.75
    \fill (h\i) circle (0.1);
    \fill[red] (hi\i) circle (0.1);
  }
  \draw (h0) node[anchor=west, node font=\tiny] {$T(1,1,1)$}
     -- (h1) node[anchor=south west, node font=\tiny] {$T(-1,1,1)$} 
     -- (h2) node[anchor=south east, node font=\tiny] {$T(-1,-1,1)$}
     -- (h3) node[anchor=east, node font=\tiny] {$T(-1,-1,-1)$}
     -- (h4) node[anchor=north east, node font=\tiny] {$T(1,-1,-1)$}
     -- (h5) node[anchor=north west, node font=\tiny] {$T(1,1,-1)$}
     -- (h0);

  % m = 8
  \draw (8,3) node[] {$m=8$};
  \foreach \i in {0,1,...,7} {
    % Vertex of 8-gon
    \coordinate (o\i) at ({8 + 1.5*cos(45*\i)}, {1.5*sin(45*\i)});  % radius = 1.5
    % Internal point (in red)
    \coordinate (oi\i) at ({8 + 0.574*cos(22.5+45*\i)}, {0.5*sin(22.5+45*\i)});  % radius = 1.5/sqrt(4+2*sqrt(2))
    
    \fill (o\i) circle (0.1);
    \fill[red] (oi\i) circle (0.1);
  }

  \draw (o0) node[anchor=west, node font=\tiny] {$T(1,1,1,1)$}
     -- (o1) node[anchor=south west, node font=\tiny] {$T(-1,1,1,1)$} 
     -- (o2) node[anchor=south, node font=\tiny] {$T(-1,-1,1,1)$}
     -- (o3) node[anchor=south east, node font=\tiny] {$T(-1,-1,-1,1)$}
     -- (o4) node[anchor=east, node font=\tiny] {$T(-1,-1,1,-1)$}
     -- (o5) node[anchor=north east, node font=\tiny] {$T(1,-1,-1,-1)$}
     -- (o6) node[anchor=north, node font=\tiny] {$T(1,1,-1,-1)$}
     -- (o7) node[anchor=north west, node font=\tiny] {$T(1,1,1,-1)$}
     -- (o0);
\end{tikzpicture}
    \caption{The zonotope reformulation for $m=6$ and $m=8$. The black dots are precisely $B_m$, which are the projection of the points $\epsilon$ satisfying the inequality constraint. The red dots are the projections of $e_k$ and $-e_k$ for $k=1,\dots \frac{m}{2}$}
    \label{fig:zonotope}
\end{figure}

An illustration of Lemma \ref{lem:zonotope_map} is shown in Figure \ref{fig:zonotope} in the cases $m=6$ and $m=8$.

We are now ready to prove Theorem \ref{thm:zonotope-reformulation} which shows how to reformulate \eqref{eq:gen-QUBO} using a binary quadratic optimization problem with $mn/2$ variables.

\begin{proof}[Proof of Theorem \ref{thm:zonotope-reformulation}]
    Let $P_m=\{(\cos(2k\pi /m),\sin(2k\pi/m):~0\le k\le m-1\}\subseteq \RR^2$ be the set of vertices of the regular $m$-gon on the real plane. We first rewrite the problem $\max_{z \in B_m^n} z^* Q z$ as a real problem on $P_m^n=\{(p_1,\ldots,p_n):~p_i\in P_m\}$: let $Q_R$ be the $2n\times 2n$ block matrix
    \[
    Q_R=\begin{bmatrix}
        Q_R(1,1) & \ldots & Q_R(1,n)\\
        \vdots & \ddots & \vdots\\
        Q_R(n,1) & \ldots & Q_R(n,n)
    \end{bmatrix},\quad  Q_R(i,j)=\begin{bmatrix}
        \text{Re}(Q_{ij}) & -\text{Im}(Q_{ij})\\
        \text{Im}(Q_{ij}) & \text{Re}(Q_{ij})
    \end{bmatrix}
    \]
    Then we have $\max_{z \in B_m^n} z^* Q z=\max_{x\in P_m^n}x^T Q_R x$. Note $Q_R$ is symmetric since $Q$ is Hermitian. 

    Let $T:\RR^{m/2}\to\CC$ be the map defined in Lemma \ref{lem:zonotope_map}. Let $\tilde{T}:\RR^{m/2}\to\RR^2$ be defined by $\tilde{T}(x)=(\text{Re}(T(x)),\text{Im}(T(x)))$, which maps $S=\{\epsilon\in\{-1,1\}^{m/2}:~\sum_{j} |\epsilon_j - \epsilon_{j+1}| \leq 2\}$ bijectively to $P_m$. In matrix form, $\tilde{T}$ is a $2\times (m/2)$ matrix given by
    \[
    A=\sin(\pi/m)\begin{bmatrix}
        \sin(\pi/m) & \sin(3\pi/m) & \ldots & \sin((m-1)\pi/m)\\
        -\cos(\pi/m) & -\cos(3\pi/m) & \ldots & -\cos((m-1)\pi/m)
    \end{bmatrix}.
    \]
    Then $I_n\otimes A$ is the matrix that maps $S^n\subseteq \{-1,1\}^{nm/2}$ bijectively to $P_m^n$, and the proof is finished. 
\end{proof}

\section{Numerical experiments}\label{sec:numerics}

In this section, we provide some numerical results to compare the performance of our zonotope-based reformulation of Theorem \ref{thm:zonotope-reformulation} vs. the standard reformulation \eqref{eq:basic-reformulation}. We solve the resultant integer programs using Gurobi \cite{gurobi} on a standard MacBook Air (M2, 2022).

\subsection{The MIMO detection problem}
\newcommand{\Symbols}{\mathsf{S}}
We consider the following problem, known as the \emph{multiple-input multiple-output (MIMO) detection problem}. Given a vector $r \in \CC^d$ defined by 
\begin{equation*}
    r : = Hx^*+\sigma v
\end{equation*}
where $H \in \CC^{d\times n}$ is the channel matrix, $v \in \CC^d$ is an additive noise vector and $\sigma \in \RR$ is the noise magnitude, our problem is to estimate the symbol signal $x^* \in \Symbols^n$, where $\Symbols$ is the set of ``symbols''.
In general, the problem can be set up with different choices for the symbol set $\Symbols$, but we only consider the relevant choice known as the $m$-ary phase shift keying ($m$-PSK), which corresponds precisely to the case where $\Symbols = B_m$ is the set of $m$'th roots of unity. 
The most common formulation for this problem, which we will work with, is based on the maximum-likelihood estimator, and is given by the following optimization problem: 
\begin{equation}\label{eq:MIMO-MLE}
    \begin{aligned}
        \text{min } & \|Hx-r\|^2\\
        \text{s.t. }&  x \in B^n_m. \\
    \end{aligned}
\end{equation}
This can be written in the form \eqref{eq:gen-QUBO} (the optimal value will be negated) with the following objective matrix:
\begin{equation}
    Q = \begin{bmatrix}
        -r^*r & r^*H \\
        H^*r & -H^*H 
    \end{bmatrix} \in \CC^{(n+1)\times (n+1)}.
\end{equation}

\begin{figure}[ht]
\centering
\begin{tikzpicture}
\begin{semilogxaxis}[enlargelimits=false,xlabel=SNR,ylabel=Running time (s),width=7cm]
\addplot+[color=blue,mark=square] table[x=snr,y=bravg]{mimo-m=4-1000.csv};
\addlegendentry{Basic reformulation}
\addplot+[color=red,mark=square] table[x=snr,y=zravg]{mimo-m=4-1000.csv};
\addlegendentry{Zonotope reformulation}
\end{semilogxaxis}
\end{tikzpicture}
\begin{tikzpicture}
\begin{semilogxaxis}[enlargelimits=false,xlabel=SNR,ylabel=Running time (s),width=7cm]
\addplot+[color=blue,mark=square] table[x=snr,y=bravg]{mimo-m=8-1000.csv};
\addlegendentry{Basic reformulation}
\addplot+[color=red,mark=square] table[x=snr,y=zravg]{mimo-m=8-1000.csv};
\addlegendentry{Zonotope reformulation}
\end{semilogxaxis}
\end{tikzpicture}
\caption{Average running time of zonotope reformulation vs. basic reformulation on a MIMO detection problem with $(d,n) = (30,20)$ and $m=4$ (left) and $m=8$ (right).}
\label{fig:MIMO-runtime}
\end{figure}
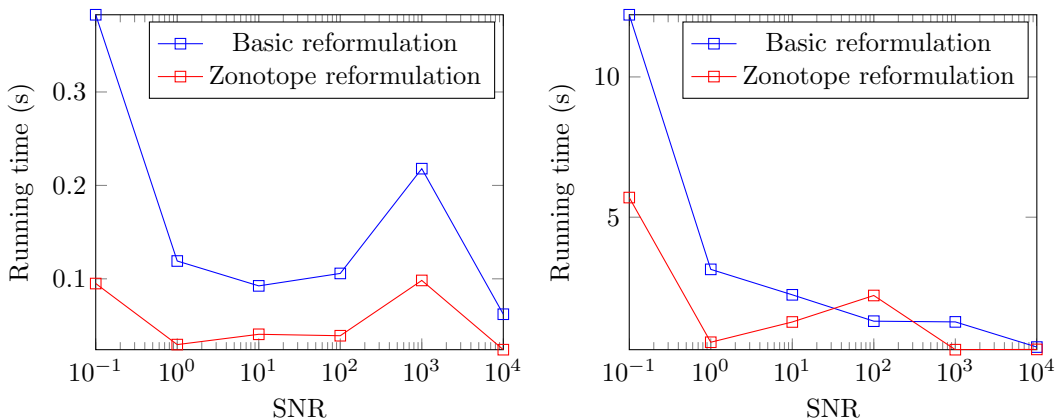

For the numerical experiments considered in this section, we consider random instances like in \cite{jiang2021tightness} where the entries of $H$ and $v$ are independent circular Gaussian variables with zero mean and unit variance. That is, the problem data will be given by $H_{ij} \sim \mathcal{CN}(0,1)$ and $v_i \sim \mathcal{CN}(0,1)$ for $i =1,\dots,d$ and $j=1,\dots,n$, where $\mathcal{CN}(0,1)$ is the standard complex Gaussian distribution. We generate the signal $x^*$ by drawing each entry $x_i^*$ for $i=1,\dots,n$ independently and uniformly from $B_m$. The signal-to-noise ratio is then defined as~\cite[Section 5]{jiang2021tightness}
\begin{equation}
    \text{SNR} = \frac{\mathbb{E}[\|H x^*\|_2^2]}{n\mathbb{E}[\|\sigma v\|_2^2]} = \frac{1}{\sigma^2}.
\end{equation}
We fix $(d,n) = (30,20)$, and consider two choices of $m$: $m =4$ and $m=8$. For each choice of $m$, we generate $1000$ instances for each SNR $=10^p$ with $p \in \{-1,0,1,2,3,4\}$, and compare the average runtime for the basic and zonotope reformulations. Both reformulations are solved using Gurobi \cite{gurobi}. The results are summarized in Figure \ref{fig:MIMO-runtime}.

We can see from Figure \ref{fig:MIMO-runtime} that the zonotope reformulation does provide an advantage over the basic reformulation in terms of solving time. In addition, the figure also demonstrates the negative correlation between the difficulty of the problem and the SNR, which was discussed in \cite{lu2019tightness,lu2020enhanced} in the context of the quality of the SDP relaxation. We expect this same phenomenon to impact the solving time of discrete optimization solvers like the one we are using, since they would depend on closing the gap between an upper bound obtained from a relaxation and a lower bound obtained from a feasible solution.

\subsection{The Potts model}

\begin{figure}[h!]
\pgfplotsset{footnotesize}
\begin{center}
\begin{tikzpicture}
\begin{semilogyaxis}[
legend columns=-1,
legend entries={Basic reformulation;, Zonotope reformulation},
legend to name=named,
title={$m = 4$}, enlargelimits=false,xlabel=$n$,ylabel=Running time (s),width=7cm, xtick = data
]
\addplot+[color=blue,mark=square] table[x=n,y=bravg]{potts-avg-runtime-m=4.csv};
\addlegendentry[color = black]{Basic reformulation}
\addplot+[color=red,mark=square] table[x=n,y=zravg]{potts-avg-runtime-m=4.csv};
\addlegendentry[color = black]{Zonotope reformulation}
\end{semilogyaxis}
\end{tikzpicture}
\begin{tikzpicture}
\begin{semilogyaxis}[
title={$m = 6$}, enlargelimits=false,xlabel=$n$,ylabel=Running time (s),width=7cm, xtick = data
]
\addplot+[color=blue,mark=square] table[x=n,y=bravg]{potts-avg-runtime-m=6.csv};
\addplot+[color=red,mark=square] table[x=n,y=zravg]{potts-avg-runtime-m=6.csv};
\end{semilogyaxis}
\end{tikzpicture}
\begin{tikzpicture}
\begin{semilogyaxis}[
title={$m = 8$}, enlargelimits=false,xlabel=$n$,ylabel=Running time (s),width=7cm, xtick = data
]
\addplot+[color=blue,mark=square] table[x=n,y=bravg]{potts-avg-runtime-m=8.csv};
\addplot+[color=red,mark=square] table[x=n,y=zravg]{potts-avg-runtime-m=8.csv};
\end{semilogyaxis}
\end{tikzpicture}
\begin{tikzpicture}
\begin{semilogyaxis}[
title={$m = 10$}, enlargelimits=false,xlabel=$n$,ylabel=Running time (s),width=7cm, xtick = data
]
\addplot+[color=blue,mark=square] table[x=n,y=bravg]{potts-avg-runtime-m=10.csv};
\addplot+[color=red,mark=square] table[x=n,y=zravg]{potts-avg-runtime-m=10.csv};
\end{semilogyaxis}
\end{tikzpicture}
\\
\ref{named}
\end{center}
\caption{Average running time of basic reformulation vs. zonotope reformulation on vector Potts model problems with the following choices of $m$ in row major order: 4, 6, 8, and 10}
\label{fig:pottsplots}
\end{figure}
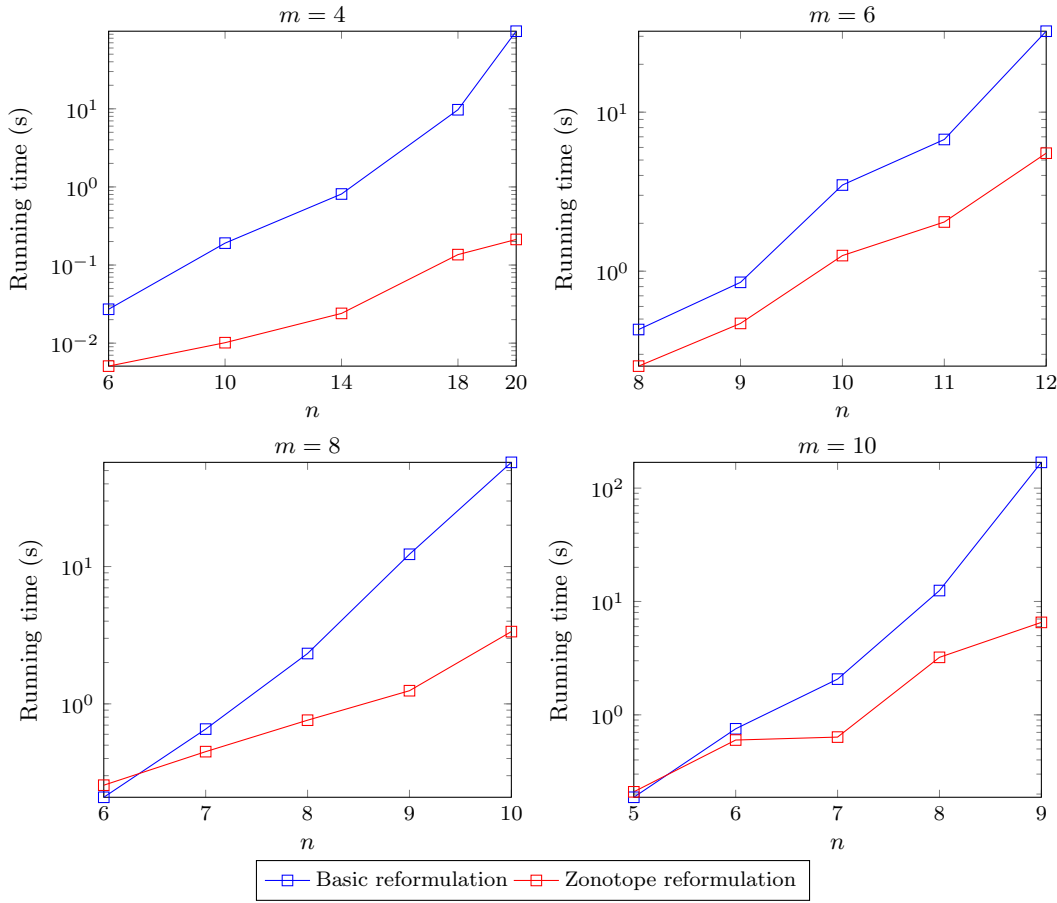

In this subsection we consider a different distribution of instances of \eqref{eq:gen-QUBO} where $Q$ is a random real symmetric matrix with entries sampled from the uniform distribution $\mathcal{U}_{[-10,10]}$. This yields instances of the ground state energy minimization problem for the vector Potts model \cite{pottsmodel} with random coupling coefficients. We use this set of instances to compare the runtime between the basic reformulation \eqref{eq:basic-reformulation} and the zonotope reformulation \eqref{eq:zonotope-reformulation}. We consider four sets of instances: $m = 4$ with $ n \in \{6,10,14,18,20\}$, $m = 6$ with $n \in \{8,9,10,11,12\}$, $m = 8$ with $n \in \{6,7,8,9,10\}$, and $m = 10$ with $n \in \{5,6,7,8,9\}$. For each set, we generate 15 instances for each choice of $n$, and record the average runtime as a function of $n$. The results are summarized in Figure \ref{fig:pottsplots}. 

The scale of the y-axis in the figures is logarithmic. Figure \ref{fig:pottsplots} demonstrates the advantage of using the zonotope reformulation instead of the basic one. We may also note that the in the case $m = 4$, this advantage is more pronounced, since in this case the zonotope reformulation is a pure QUBO (i.e., unconstrained).

\section{Conclusions}\label{sec:conclusions}

In this paper, we studied the moment-SOS hierarchy in the context of the generalized QUBO problem \eqref{eq:gen-QUBO}, also known as discrete complex quadratic optimization. We studied the sum-of-squares hierarchy for every $m$ and showed that it converges at level $\lceil n/2\rceil$, generalizing the result for $m=2$ \cite{fawzi2016sparse}. 

In addition, we have constructed, for even $m$, a new real binary reformulation of \eqref{eq:gen-QUBO} based on zonotopes, where the number of binary variables is $\frac{mn}{2}$ instead of the $mn$ variables one would get from the basic real binary reformulation. The two reformulations were tested and compared on MIMO detection and vector Potts model instances. The zonotope reformulation exhibits significant advantage in runtime.
%, especially for $m=4$ where the reformulation is unconstrained, i.e., pure QUBO. \\

In terms of future work, there are two main directions that we think form a natural extension to the work in this paper.

The first question concerns Theorem \ref{thm:main}: our theorem works for the general class of polynomials $f$ whose terms only depend on two variables each. An interesting question is to know whether one can improve the result when the monomials are all of the form $z_i^{-1} z_j$, by using an alternative grading of monomials which has fewer terms at each level. A related question is to know whether the level $\lfloor n/2 \rfloor+1$ is tight or can be improved by constructing explicit nonnegative polynomials $f$ when $m \geq 3$ that require high-degree sum-of-squares proofs.

The second direction would be to provide an alternative to the basic real binary reformulation \eqref{eq:basic-reformulation} for the case where $m$ is odd. The zonotope reformulation we have provided in this work is only valid for even $m$, and so a natural next step would be to find a unified reformulation recipe that would work for all $m \geq 2$.

\paragraph{Acknowledgments} AA acknowledges funding from the Kingdom of Saudi Arabia's Ministry of Education under the research and development scholarship track. HF and SS acknowledge funding from
UK Research and Innovation (UKRI) under the UK government’s Horizon Europe funding guarantee
EP/X032051/1.

\bibliography{bibs}
\bibliographystyle{siam}

\end{document}